\documentclass[10pt]{amsart}
\usepackage{amsmath, amsthm, amsfonts, amssymb, graphicx, hyperref, bm,verbatim,mathrsfs,siunitx,natbib,tikz-cd}
\usepackage{stmaryrd,enumerate}
\theoremstyle{plain}
\newtheorem{thm}{Theorem}[section]
\newtheorem{lem}[thm]{Lemma}
\newtheorem{prop}[thm]{Proposition}
\newtheorem{cor}[thm]{Corollary}

\numberwithin{sublem}{thm} 
\numberwithin{equation}{section}

\newcommand{\e}{\varepsilon}
\begin{document}
\title{Simultaneously Small Fractional Parts of Polynomials}
\author{Cheuk Fung (Joshua) Lau}
\address{University College, Oxford OX1 4BH, UK}
\email{joshua.lau@univ.ox.ac.uk}
\begin{abstract}
Let $f_1,\dots,f_k\in\mathbb{R}[X]$ be polynomials of degree at most $d$ with $f_1(0)=\dots=f_k(0)=0$. We show that there is an $n<x$ such that $\|f_i(n)\|\ll x^{-1/10.5kd(d-1)+o(1)}$ for all $1\le i\le k$. This improves on an earlier result of Maynard, who obtained the same exponent dependency on $k$ but not on $d$.
\end{abstract}
\maketitle
%
%
%
%
%
%
%
\section{Introduction}
The question of how small the fractional part $\lVert f(n) \rVert$ of a polynomial $f(n)$ (with $f(0)=0$) can be made has been investigated since the early 20th century, with the current record being
\begin{equation}
\min_{n \leq x} \ \lVert f(n) \rVert \ll_d x^{-1/2d(d-1)+o(1)}, \label{eq:1DimBest}
\end{equation}
which is due to \cite{baker2016small} (for $d \geq 8$). Here
$\|\cdot\|$ denotes the distance to the nearest integer. We impose the condition $f(0)=0$ so as to avoid examples like $f(n)=n+1/2$, which clearly doesn't attain arbitrarily small fractional parts. Note that the bound depends only on $x$ and $d$, and is otherwise completely uniform over all such polynomials $f$. The exponent $1/2d(d-1)$ is based on proving Vinogradov's Mean Value Theorem for $k\geq 4$ by \citet{bourgain2016proof}. This bound is not expected to be tight, and it is conjectured \citep{baker1986diophantine} that this should be improvable to $1+o(1)$.

One of the first results of this kind is from \citet{heilbronn1948distribution}, who obtained for any real $\alpha$, $\e>0$ and sufficiently large $x>0$, 
$$\min_{n \leq x} \lVert \alpha n^2 \rVert<x^{-1/2+\e}.$$
For the problem of a general polynomial $f \in \mathbb{R}[X]$ vanishing at the origin, \citet{selberg1955vinogradov} obtained
$$\min_{n \leq x} \lVert f(n) \rVert \ll_d x^{-1/(8+o(1)) d^2 \log d}.$$
\citet{wooley2012vinogradov} improved the exponent to $-1/4d(d-1)$ for $d \geq 2$, and as noted \citet{baker2016small} improved it to $-1/2d(d-1)$ for $d \geq 8$.

We can ask the same question, but instead for $k$ polynomials $f_1,\dots,f_k\in\mathbb{R}[X]$ of degree at most $d$ with $f_1(0)=\dots =f_k(0)=0$, the current record due to \cite{maynard2021simultaneous} is for some $c_d>0$ depending only on $d$,
\begin{equation}
\min_{n\le x}\max_{i\le k}\|f_i(n)\|\ll_{k,d} \frac{1}{x^{c_d/k}}.\label{eq:ManyPolysMaynard}
\end{equation}
This refines the work of \cite{baker1980fractional}, where he obtained a hybrid bound
\begin{equation}
\min_{n\le x}\max_{i\le k}\|f_i(n)\|\ll_{k,d} \frac{1}{x^{1/(k^2+k c_d)+o(1)}}.\label{eq:ManyPolysBaker}
\end{equation} A heuristic based on choosing the coefficients of $f_1,\dots,f_k$ uniformly at random shows that one could only expect
\begin{equation}
\min_{n\le x}\max_{i\le k}\|f_i(n)\|\ll \frac{1}{x^{1/k}}.\label{eq:Optimal}
\end{equation}

Our main result is to establish a bound for \eqref{eq:ManyPolysMaynard} with an exponent of the form $1/10.5kd(d-1)+o(1)$. By comparing this with \eqref{eq:Optimal} and \eqref{eq:1DimBest}, we see that this bound is best in the $k$ and $d$ aspect, apart from hybrid bounds and the constant 10.5. More precisely, our main result is the following.
\begin{thm}\label{thm:MainTheorem}
Let $k, d \in \mathbb{Z}^+,\, \e>0,$ and $M \in \mathbb{R}^+$ satisfy
$$M \geq \max \left\{4, \frac{1}{2}+\frac{\log \e^{-1}}{2 \log 2} \right\}.$$
There is a constant $C_{d, k}>2$ depending only on $d$, $k$ (and $\e,M$) such that the following holds.

Let $f_1, \ldots, f_k \in \mathbb{R}[X]$ be polynomials of degree at most $d$ such that $f_1(0)=\cdots=$ $f_k(0)=0$. Let $\e_1, \ldots, \e_k \in(0,1 / 100]$, and put $\Delta=\prod_{i=1}^k \e_i$. Define $$c_2=10.5+9k/(2k)^M+\e.$$

If $\Delta^{-1} \leq x^{1 / c_2d(d-1)}$ and $x>C_{d, k}$, then there is a positive integer $n<x$ such that
$$
\left\|f_i(n)\right\| \leq \e_i \quad \text { for all } i \in\{1, \ldots, k\} \text {. }
$$
\end{thm}
%
%
%
%
By taking $\e_1=\cdots=\e_k=x^{-1/c_2d(d-1)k}$, and for any $\e>0$ letting $M$ large such that $2^M \geq \e^{-1}$, we arrive at the following:
\begin{cor}
    Let $f_1,\ldots,f_k \in \mathbb{R}[X]$ be polynomials of degree at most $d$ such that $f_1(0)=\cdots=f_k(0)=0$. Then there is a positive integer $n<x$ such that for all $\e>0$,
    $$\lVert f_i(n) \rVert \ll_{d,k,\e} x^{-1/(10.5+\e)kd(d-1)}$$
    for all $i=1,2,\ldots,k$.
\end{cor}
%
%
%
%
The proof of Theorem \ref{thm:MainTheorem} is obtained by following the argument of \cite{maynard2021simultaneous}, but making more careful choices of constants in the proofs. If we follow strictly the argument based on Weyl differencing in \cite{maynard2021simultaneous}, $c_d$ would be of the form $1/d^d$. To improve on the $d$-quantification, one can use arguments based on the Vinogradov Mean Value Theorem such as \cite{wooley2012vinogradov}, which gives $c_d$ to be of the form $1/2^d$. The main improvement made is found in Lemma \ref{alotorsamedenom}, which allows us to take $c_d$ to be of the form $1/d^2$.
%
%
%
%

%
%
%
%
%
%
%
%
\section{Outline}
In this section we outline the main improvement in this paper. The main method follows closely the argument of \citet{maynard2021simultaneous}, where in Section \ref{smallfractionalpartsormanyrelations} we prove a version of \citet[][Proposition 5.1]{maynard2021simultaneous}, in Section \ref{manyrelationshavesamedenominator} we prove a version of \citet[][Proposition 5.2]{maynard2021simultaneous}, and in Section \ref{densityincrement} we prove a version of \citet[][Proposition 5.3]{maynard2021simultaneous}. Sections \ref{smallfractionalpartsormanyrelations} and \ref{densityincrement} are very similar to the corresponding sections in \citet{maynard2021simultaneous}. To emphasize on the improvement made in Section \ref{manyrelationshavesamedenominator}, we first review the arguments in the proof of \cite[][Lemma 7.3]{maynard2021simultaneous}.\\

Given $f_i(X)=\sum_{j=1}^d f_{i,j} X^j \in \mathbb{R}[X]$, Fourier analysis shows that, given any intervals $I_1,\ldots,I_k$ of length $\e>0$ (for some small $\e>0$) and $x>\e^{-kC_d}$, either there is an $n<x$ such that the vector of fractional parts $\mathbf{v}(n)=(\lVert f_1(n) \rVert,\ldots,\lVert f_k(n) \rVert)$ lie in $I_1 \times \cdots \times I_k$, or there are at least $Q^{1/cd(d-1)}$ integer tuples $(h_1,\ldots,h_k)$ with $h_i \leq \e^{-1-o(1)}$ such that
\begin{equation} \label{outlinerationalapprox}
    h_1f_{1,j}+\cdots+h_kf_{k,j} \approx \frac{a_j}{q_j} \text{ for all } 1 \leq j \leq d,
\end{equation}
for some $a_j \in \mathbb{Z}$ and $Q \leq \e^{-ckd(d-1)}$ with $q_j<Q^{1/cd}$. \citet{maynard2021simultaneous} proved that, given any $\delta>0$ small, $r \in \mathbb{Z}^+$ and $\mathcal{S} \subseteq \mathbb{Z} \times \mathbb{Z} \times \mathbb{Z}^k$ set of triples $(a,q,\mathbf{h})$ with $\gcd(a,q)=1$ and $q \leq Q^{1/cd}$ such that $\# \mathcal{S} \geq Q^{\delta}$, then one of the following holds:
\begin{enumerate}
    \item At least $\# \mathcal{S}^{1/2}$ of triples $(a,q,\mathbf{h})$ have the same $q=q_0$.
    \item There are at least $\# \mathcal{S}^{r/5}$ distinct values of $a_1/q_1+\cdots+a_r/q_r$, amongst all $(a_i,q_i,\mathbf{h}^{(i)}) \in \mathcal{S}$.
\end{enumerate}
Indeed, if $(1)$ is false, then one can find an integer $m_0$ and a subset of the $q$s such that they are all multiples of $m_0$, and not too many of them are multiples of $m_0 \ell$ for $\ell>1$. Letting $b=q/m_0$ for such $m_0$, one can prove that the size of the set of all tuples $(b_1,\ldots,b_r)$ with pairwise small gcd is not small. Focusing on these tuples, one can prove that there are many distinct values of $a_1/m_0b_1+\cdots+a_r/m_0b_r$.\\

Therefore, if we are in the case of (\ref{outlinerationalapprox}), we can focus on the case $j=1$ to prove there are many possible values of $\sum a_1/q_1$, hence there are at least $Q^{1/2cd(d-1)}$ triples $(\mathbf{a},\mathbf{q},\mathbf{h}) \in \mathbb{Z}^d \times \mathbb{Z}^d \times \mathbb{Z}^k$ satisfying (\ref{outlinerationalapprox}), with $q_1$ being identical amongst all triples. Doing this procedure $d$ times, there are at least $Q^{1/2^dcd(d-1)}$ triples $(\mathbf{a},\mathbf{q},\mathbf{h})$ satisfying (\ref{outlinerationalapprox}) with $\mathbf{q}$ identical amongst all triples. This eventually leads to an exponential dependence on $d$ in the final exponent.\\

To get around this issue, we consider a multi-dimensional argument, instead of iterating a $1$-dimensional argument $d$ times. We prove that given $\mathcal{S} \subseteq \mathbb{Z}^d \times \mathbb{Z}^d \times \mathbb{Z}^k$ set of triples $(\mathbf{a},\mathbf{q},\mathbf{h})$ with $\gcd(a_i,q_i)=1$ and $q_i \leq Q^{1/cd}$ (for $i=1,2,\ldots,d$) such that $\# \mathcal{S} \geq Q^{\delta}$, then one of the following holds:
\begin{enumerate}
    \item At least $\# \mathcal{S}^{1/2}$ of triples $(\mathbf{a},\mathbf{q},\mathbf{h})$ have the same $\mathbf{q}=\mathbf{q}^{(0)}$.
    \item There are at least $\# \mathcal{S}^{r/5}$ distinct values of $$\left( \frac{a_1^{(1)}}{q_1^{(1)}}+\cdots+\frac{a_1^{(r)}}{q_1^{(r)}},\ldots, \frac{a_d^{(1)}}{q_d^{(1)}}+\cdots+\frac{a_d^{(r)}}{q_d^{(r)}} \right),$$ amongst all $(\mathbf{a}^{(i)},\mathbf{q}^{(i)},\mathbf{h}^{(i)}) \in \mathcal{S}$.
\end{enumerate}
Indeed, if $(1)$ is false, then one can find an integer tuple $\mathbf{m}^{(0)}$ and a subset of the $\mathbf{q}$s such that all $q_i$ are multiples of $m_i^{(0)}$, and not too many of them are multiples of $m_i^{(0)} \ell_i$, where $\bm{\ell} \neq (1,\ldots,1)$. Letting $b_i := q_i/m_i^{(0)}$ for such $\mathbf{m}^{(0)}$, one can prove that the size of the set of all tuples $(\mathbf{b}^{(1)},\ldots,\mathbf{b}^{(r)})$ with pairwise small gcd (at all entries) is not small. Focusing on these tuples, we can prove that there are many distinct vectors of the aforementioned form.\\

Therefore, if we have (\ref{outlinerationalapprox}), we can prove there are many possibilities of $(\sum a_1/q_1,\ldots,\sum a_d/q_d)$, hence there are at least $Q^{1/2cd(d-1)}$ triples $(\mathbf{a},\mathbf{q},\mathbf{h}) \in \mathbb{Z}^d \times \mathbb{Z}^d \times \mathbb{Z}^k$ satisfying (\ref{outlinerationalapprox}), with $\mathbf{q}$ identical amongst all triples. This multi-dimensional argument avoids the $2^d$ cost in the exponent, which eventually allows a quadratic dependency on $d$ in the final exponent.

\section{Acknowledgements}
%
%
%
%
We would like to thank James Maynard for suggesting this question, and also the guidance we received throughout the writing of this paper.
%
%
%
%
\section{Notation}
%
%
%
%
Throughout the paper we assume that we have polynomials $f_1,\dots,f_k\in\mathbb{R}[X]$ of degree at most $d$ with $f_1(0)=\dots=f_k(0)=0$. We let these polynomials be given by $f_i(X)=\sum_{j=1}^d f_{i,j}X^j$. Furthermore, we have reals $\e_1,\dots,\e_k\in(0,1/100]$, and we put $\Delta:=\prod_{i=1}^k\e_i$.

We often use Vinogradov's notation, where $X \ll Y$ denotes $|X| \leq CY$ for some constant $C$. If the implied constant depends on $\e$ say, we write $X \ll_\e Y$. We also use Landau's Big O notation, where $X=O(Y)$ if $X \ll Y$.
%
%
%
%
\section{Small Fractional Parts or Many Relations} \label{smallfractionalpartsormanyrelations}
%
%
%
%
Throughout this paper, we follow closely the arguments of \cite{maynard2021simultaneous}. For the sake of clarity, we include the proofs for statements that are not identical to the analogous ones in \cite{maynard2021simultaneous}, whereas for identical results we just state without proof.
\begin{lem}
\label{equi10}
     Let $f(X)=\sum_{i=1}^d f_iX^i \in \mathbb{R}[X]$ be a polynomial of degree $d \geq 2$ with $f(0)=0$ and $\e>0$. If there is some $Q \in [2,x^{1/2d(d-1)-\e}]$ such that $$\left| \sum_{n<x} e(f(n))\right| \geq \frac{x}{Q},$$
    then there is an integer $q<x^\e Q^{d}$ and integers $a_1,\ldots,a_d$ such that
    $$f_j=\frac{a_j}{q}+O \left( \frac{Q^{d}}{x^{j-\e}} \right)$$ for $j=1,2,\ldots,d$. The implied constant only depends on $d$ and $\e$.
\end{lem}
\begin{proof}
    Follows immediately from \citet[][Theorem 4]{baker2016small}.
\end{proof}
\begin{lem}
\label{6.3}
    Let $\e>0,M \geq 4$, $f_1,\ldots,f_k \in \mathbb{R}[X]$ be real valued functions, and $\e_1,\ldots,\e_k \in (0,1/2]$ be real numbers, with $\Delta := \prod_{i=1}^k \e_i$. Suppose $x$ is sufficiently large in terms of $\e$. Then at least one of the following holds:
    \begin{enumerate}
        \item We have
        $$\# \{ n \leq x:\lVert f_i(n) \rVert \leq \e_i \ \forall i \} >0.$$
        \item There is a quantity $Q \geq 2$ of the form $Q=2^{j}$ such that there are at least $Q^{1/(1+\e)}$ distinct values of $\mathbf{h} \in \mathbb{Z}^k \backslash \{\mathbf{0}\}$ with $|h_i|<\e_i^{-1} \Delta^{-1/(2k)^M}$ such that
        $$\frac{x}{Q} \leq \left| \sum_{n \leq x} e \left( \sum_{i=1}^k h_if_i(n) \right) \right| \leq \frac{2x}{Q}.$$
    \end{enumerate}
\end{lem}
\begin{proof}
    The proof is nearly identical to that of \citet[][Lemma 6.3]{maynard2021simultaneous}. We fix a smooth function $\phi:\mathbb{R} \to [0,1]$ with $\phi(t)$ supported on $|t|<1$, which is $1$ on $|t|<1/2$. Let
    $$\Phi_i(t)=\sum_{m \in \mathbb{Z}} \phi \left( \frac{t+m}{\e_i} \right),$$
    which is clearly 1-periodic, smooth, and supported on $\lVert t \rVert<\e_i$. By Poisson summation,
    $$\Phi_i(f_i(t))=\e_i \sum_{h \in \mathbb{Z}} \hat{\phi}(\e_ih)e(hf_i(t)).$$
    Since $\phi$ is fixed and smooth, $\phi^{(j)}(t) \ll_j 1$, so $|\hat{\phi}(u)| \ll_j u^{-j}$ for $j \geq 0$. Thus we see that the terms with $|h| \geq \e_i^{-1} \Delta^{-1/(2k)^M}$ contribute $O(\Delta^{100})$, and so
    $$\Phi_i(f(t))=\e_i \sum_{|h| \leq \e_i^{-1} \Delta^{-1/(2k)^M}} \hat{\phi}(\e_ih)e(hf_i(t))+O(\Delta^{100}).$$
    Thus we find that
    \begin{align*}
        \# \{n \leq x:\lVert f_i(n) \rVert \leq \e_i \forall i\} &\geq \sum_{n \leq x} \prod_{i=1}^k \Phi_i(f_i(n))\\
        &= \Delta \sum_{\substack{h_1,\ldots,h_k\\ |h_i|<\e_i^{-1} \Delta^{-1/(2k)^M}}} \left( \prod_{i=1}^k \hat{\phi}(\e_ih_i) \right) \sum_{n \leq x} e \left( \sum_{i=1}^k h_if_i(n) \right)+O(x \Delta^{99})\\
        &= x \Delta \hat{\phi}(0)^k+O \left( \Delta \sum_{\substack{\mathbf{h} \in \mathbb{Z}^k \setminus \{0\}\\ |h_i|<\e_i^{-1} \Delta^{-1/(2k)^M}}} \left| \sum_{n \leq x} e \left( \sum_{i=1}^k h_if_i(n) \right) \right| \right) + O(x \Delta^{99}).
    \end{align*}
    For $\Delta$ sufficiently small, we see that $\Delta \hat{\phi}(0)^k+O(\Delta^{99}) \gg \Delta$, and so either
    $$\# \{n \leq x: \lVert f_i(n) \rVert \leq \e_i \forall i\} >0,$$
    or $$\sum_{\substack{\mathbf{h} \in \mathbb{Z}^k \setminus \{0\}\\ |h_i|<\e_i^{-1} \Delta^{-1/(2k)^M}}} \left| \sum_{n \leq x} e \left( \sum_{i=1}^k h_if_i(n) \right) \right| \gg x.$$
    In the latter case, by pigeonhole principle (or more specifically \citet[][Lemma 2.1]{bloom2020breaking}), there is some $Q=2^{j}$ such that there are at least $Q^{1-\e}$ choices of $\mathbf{h}$ in the outer summation satisfying
    $$\frac{x}{Q} \leq \left| \sum_{n \leq x} e \left( \sum_{i=1}^k h_if_i(n) \right) \right| \leq \frac{2x}{Q}.$$
\end{proof}
 Using this, we can prove the analogous version of \cite[][Proposition 5.1]{maynard2021simultaneous}.
\begin{prop}
\label{5.1}
    Let $f_1,\ldots,f_k \in \mathbb{R}[X]$ be polynomials of degree at most $d$ such that $f_1(0)=\cdots=f_k(0)=0$. Put $f_i(X)=\sum_{j=1}^d f_{i,j} X^j$. Let $\e_1,\ldots,\e_k \in (0,1/100]$, and put $\Delta=\prod_{i=1}^k \e_i$. Let $\e>0,M \geq 4$. Define constants
    $$c_0=1+\e, \quad c_1=1+k/(2k)^{M}+10\e, \quad c_2=2+1/(2k)^{M-1}+20\e.$$

     Then provided $\Delta^{-c_2d(d-1)}<x$, at least one of the following holds:
    \begin{enumerate}
        \item We have
        $$\# \{n \leq x:\lVert f_i(n) \rVert <\e_i \ \forall i \} >0.$$
        \item There is some $Q \leq \Delta^{-c_1d(d-1)}$ such that there are at least $Q^{1/(1+\e)c_0d(d-1)}$ triples $(\mathbf{a},\mathbf{q},\mathbf{h}) \in \mathbb{Z}^d \times \mathbb{Z}^d \times \mathbb{Z}^k$ satisfying:
        \begin{enumerate}
            \item $q_1=\cdots=q_d$ and $1 \leq q_j \leq x^\e Q^{1/c_0(d-1)}$ for $1 \leq j \leq d$.
            \item $h_i \ll \e_i^{-1} \Delta^{-1/(2k)^M}$ for $1 \leq i \leq k$.
            \item For each $j=1,2,\ldots,d$, we have
            $$\sum_{i=1}^k h_{i} f_{i,j}=\frac{a_j}{q_j}+O \left( \frac{Q^{1/d}}{x^{j-\e}} \right).$$
        \end{enumerate}
    \end{enumerate}
 All implied constants depend only on $d,k,\e$ and $M$..
\end{prop}
\begin{proof}
    The proof is very similar to that of \citet[][Proposition 5.1]{maynard2021simultaneous}. Suppose condition (1) does not hold. Then, by Lemma \ref{6.3}, there is $Q_1 \geq 2$ such that there are at least $Q_1^{1/(1+\e)}$ distinct values of $\mathbf{h} \in \mathbb{Z}^k \backslash \{ \mathbf{0} \}$ with $|h_i|<\e_i^{-1} \Delta^{-1/(2k)^M}$ such that
    $$\left| \sum_{n \leq x} e \left( \sum_{i=1}^k h_i f_i(n) \right) \right| \geq \frac{x}{Q_1}.$$
    Letting $Q=Q_1^{c_0d(d-1)}$, there are at least $Q^{1/(1+\e)c_0d(d-1)}$ distinct values of $\mathbf{h} \in \mathbb{Z}^k \backslash \{ \mathbf{0} \}$ with $|h_i|<\e_i^{-1} \Delta^{-1/(2k)^M}$ such that
    $$\left| \sum_{n \leq x} e \left( \sum_{i=1}^k h_i f_i(n) \right) \right| \geq \frac{x}{Q^{1/c_0d(d-1)}}$$
    Note that $Q_1 \leq \Delta^{(-1-k/(2k)^M)(1+\e)}$, so we have $Q_1^{c_0d(d-1)}=Q \leq \Delta^{-c_1d(d-1)}$ since $c_1 \geq c_0(1+k/(2k)^M)(1+\e)$, and
    $$Q_1=Q^{1/c_0d(d-1)} \leq \Delta^{-(1+k/(2k)^M)(1+\e)} \leq \Delta^{-(c_2+2\e)d(d-1)/2d(d-1)} < x^{1/2d(d-1)-\e},$$
    since $c_2\geq 2+1/(2k)^{M-1}+20\e$.
    Therefore applying Lemma \ref{equi10} to $Q_1$, there are integers $q_1,\ldots,q_d< x^\e Q_1^{d}= x^\e Q^{1/c_0(d-1)}$ and integers $a_1,\ldots,a_d$ such that $\gcd(a_j,q_j)=1$ for all $j=1,2,\ldots,d$ and
    $$f_j=\frac{a_j}{q_j} + O \left( \frac{Q^{1/d}}{x^{j-\e}} \right)$$
    for $j=1,2,\ldots,d$. 
\end{proof}
\section{Many Relations Have Same Denominator} \label{manyrelationshavesamedenominator}
The proof of the following lemma can be found in the proof of \cite[][Lemma 7.1]{maynard2021simultaneous}. We state it below for sake of clarity.
\begin{lem} \label{smallgcdfewdenom}
    Let $B>1,\delta>0,$ and $r \in \mathbb{Z}^+$. For $b_1,\ldots,b_r \in \mathbb{Z}$ satisfying
    $$B \leq b_i<2B, \quad \gcd(b_i,b_j)<B^{\delta^2/r^2} \text{ for all } 1 \leq i<j \leq r,$$
    define the set $\mathcal{R}(b_1,\ldots,b_r)$ consisting of $b_1',\ldots,b_r' \in \mathbb{Z}$ satisfying
    \begin{enumerate}
        \item $B \leq b_i' <2B,$ $\gcd(b_i',b_j')<B^{\delta^2/r^2} \text{ for all } 1 \leq i<j \leq r,$
        \item there exists $a_1,\ldots,a_r,a_1',\ldots,a_r' \in \mathbb{Z}$ such that
        \begin{enumerate}
            \item $\gcd(a_i,b_i),\gcd(a_i',b_i')=1$ for all $1 \leq i \leq r$,
            \item $\frac{a_1}{b_1}+\cdots+\frac{a_r}{b_r}=\frac{a_1'}{b_1'}+\cdots+\frac{a_r'}{b_r'}$.
        \end{enumerate}
    \end{enumerate}
    Then $\# \mathcal{R}(b_1,\ldots,b_r) \leq B^{5\delta^2}$.
\end{lem}
\begin{proof}
    For any choice of $a_1,\ldots,a_r$ with $\gcd(a_i,b_i)=1$, the denominator of $a_1/b_1+\cdots+a_r/b_r$ is a multiple of $p^\ell$ if $p^\ell$ divides exactly one of $b_1,\ldots,b_r$, and $p^{\ell+1}$ divides none of them. Let $\gcd(b,p^\infty)$ denote the largest power of $p$ dividing $b>1$, and $\gcd(b_i,b_j,p^\infty)$ the largest power of $p$ dividing both $b_i$ and $b_j$. Define
    $$g_p := \frac{\prod_{i=1}^r \gcd(b_i,p^\infty)}{\prod_{1 \leq i <j \leq r} \gcd(b_i,b_j,p^\infty)^2}.$$
    We see that if $p^\ell$ divides exactly one of $b_1,\ldots,b_r$, say $b_1$, and $p^{\ell+1}$ divides none of them, then $g_p \leq p^\ell$ since $\gcd(b_1,b_j,p^\infty)^2 \geq \gcd(b_j,p^\infty)$ for all $1 < j \leq r$. Similarly, if $p^\ell$ divides at least 2 of the $b_i$ but $p^{\ell+1}$ divides none of them, then $g_p \leq 1$.\\

    Therefore, for any $a_1,\ldots,a_r$ with $\gcd(a_i,b_i)=1$, the denominator of $a_1/b_1+\cdots+a_r/b_r$ is of size at least $\prod_p g_p$ by considering the prime factorisation of $b_i$. However,
    $$\prod_p g_p=\frac{\prod_{i=1}^r b_i}{\prod_{1 \leq i<j \leq r} \gcd(b_i,b_j)^2} \geq B^{r-2\delta^2}.$$
    Clearly any such denominator is of size $O(B^r)$, so there are $O(B^{2 \delta^2})$ possible denominators of $a_1/b_1+\cdots+a_r/b_r$. Given such a denominator $q>B^{r-2\delta^2}$, if it is also the denominator of $a_1'/b_1'+\cdots+a_r'/b_r'$, then $q$ divides $\prod_{i=1}^r b_i'$.\\

    As $\prod_{i=1}^r b_i' \ll B^r$, this implies there are $O(B^{2 \delta^2})$ choices of $\prod_{i=1}^r b_i'$, so $O(B^{2\delta^2+o(1)})$ choices of $b_1',\ldots,b_r'$ using the divisor bound. Therefore given $(b_1,\ldots,b_r)$ satisfying the above requirements, there are at most $B^{5 \delta^2}$ choices of $(b_1',\ldots,b_r')$ in total, so $\# \mathcal{R}(b_1,\ldots,b_r) \leq B^{5 \delta^2}$.
\end{proof}
 We prove a generalisation of \cite[][Lemma 7.1]{maynard2021simultaneous}. In \cite{maynard2021simultaneous}, Lemma 7.3 was proven with induction, which introduced extra factors of $d$ on the exponent. We directly prove a generalisation which avoids this issue.
\begin{lem}{(Expansion or Same Denominators)} \label{alotorsamedenom}
    Let $\delta,\e \in (0,1)$, $r$ a positive integer, $d \geq 2$, and $X \in \mathbb{R}_{\geq 1}$. Let $Q>0$ be sufficiently large in terms of $r, \ d$, $\delta$ and $\e$, and let $\mathcal{S} \subseteq \mathbb{Z}^d \times \mathbb{Z}^d \times \mathbb{Z}^k$ be a set of triples $(\mathbf{a},\mathbf{q},\mathbf{h})$ with $\gcd(a_i,q_i)=1$ and $q_i \leq  XQ^{1/(d-1)}$ (for $i=1,2,\ldots,d$) such that $\# \mathcal{S} \geq Q^{\delta}$.\\
    
    Then one of the following holds:
    \begin{enumerate}
        \item There is a $\mathbf{q}^{(0)}$ such that at least $\# \mathcal{S}^{1/2}$ of the triples $(\mathbf{a},\mathbf{q},\mathbf{h}) \in \mathcal{S}$ have $\mathbf{q}=\mathbf{q}^{(0)}$.
        \item The set 
        \footnotesize
        $$\mathcal{A} := \left\{ 
\left( \frac{a_1^{(1)}}{q_1^{(1)}}+\cdots+\frac{a_1^{(r)}}{q_1^{(r)}},\ldots, \frac{a_d^{(1)}}{q_d^{(1)}}+\cdots+\frac{a_d^{(r)}}{q_d^{(r)}} \right) : \exists \mathbf{h}^{(1)},\ldots,\mathbf{h}^{(r)} \in \mathbb{Z}^k \text{ s.t. } (\mathbf{a}^{(i)},\mathbf{q}^{(i)},\mathbf{h}^{(i)}) \in \mathcal{S} \text{ for } 1 \leq i \leq r \right\},$$
\normalsize (where $\mathbf{a}^{(i)}=(a_1^{(i)},\ldots,a_d^{(i)})$ and $\mathbf{q}^{(i)}=(q_1^{(i)},\ldots,q_d^{(i)}))$ has cardinality at least $X^{-\e \delta dr/40-5\e \delta^2d/r} \# \mathcal{S}^{(\frac{1}{2}-\e)r}$.
    \end{enumerate}
\end{lem}
\begin{proof}
    The proof is a generalisation of Lemma {7.1} found in \cite{maynard2021simultaneous}. Throughout the lemma we will assume that $Q$ is large enough in terms of $\e, \delta$ and $r$. We first restrict our attention to a suitable subset of the $\mathbf{q}'s$ appearing in $\mathcal{S}$. For $j_1,j_2,\ldots,j_d \in \mathbb{N}_0$, let
    $$\mathcal{B}_{j_1,\ldots,j_d}=\left\{\mathbf{q} \in \prod_{i=1}^d [2^{j_i},2^{j_i+1}): \exists (\mathbf{a},\mathbf{h}) \in \mathbb{Z}^d \times \mathbb{Z}^k \text{ with } \gcd(a_i,q_i)=1 \ \forall i \text{ and } (\mathbf{a},\mathbf{q},\mathbf{h}) \in \mathcal{S} \right\}.$$
    Clearly $\mathcal{B}_{j_1,\ldots,j_d}=\emptyset$ if $j_i>2 \log (XQ^{2/d})$ for some $i$. Note if
    $$\# \left\{ \mathbf{q}:\exists (\mathbf{a},\mathbf{h}) \in \mathbb{Z}^d \times \mathbb{Z}^k \text{ with } (\mathbf{a},\mathbf{q},\mathbf{h}) \in \mathcal{S} \right\}=\sum_{2^{j_i} \leq Q^{1/d} \, \forall i} \# \mathcal{B}_{j_1,\ldots,j_d} \leq \# \mathcal{S}^{1/2}$$
    then by pigeonhole principle, there is a $\mathbf{q} \in \mathbb{Z}^d \cap [1, XQ^{2/d}]^d$ such that there are at least $\# \mathcal{S}^{1/2}$ choices of $(\mathbf{a},\mathbf{h})$ with $(\mathbf{a},\mathbf{q},\mathbf{h}) \in \mathcal{S}$. Thus condition (1) is satisfied in this case.\\

     Therefore we may assume that
    $$\sum_{j_i} \# \mathcal{B}_{j_1,\ldots,j_d} > \# \mathcal{S}^{1/2} \geq Q^{\delta/2},$$
    and so there is some $j_1,\ldots,j_d \leq 2 \log 
    (XQ^{2/d})$ such that $\# \mathcal{B}_{j_1,\ldots,j_d}>\# \mathcal{S}^{\frac{1}{2}-\e_1}$, where $\e_1=\e/20$. Note we must have $2^{j_1+\cdots+j_d}>Q^{\delta/3}$ from the trivial bound $\# \mathcal{B}_j \leq 2^{j_1+\cdots+j_d}$.\\

     For convenience, we define the following terminology. We say $\mathbf{m} \in \mathbb{Z}^d$ divides $\mathbf{q} \in \mathbb{Z}^d$ if $m_i \mid q_i$ for all $1 \leq i \leq d$, and define $\mathbf{m} \mathbf{q} := (m_1q_1,\ldots,m_dq_d)$.\\

     If there is a $d$-tuple $\mathbf{m}$ which divides at least $\# \mathcal{B}_{j_1,\ldots,j_d}/m_1^{\e_1 \delta/2} \cdots m_d^{\e_1 \delta/2}$ elements of $\mathcal{B}_{j_1,\ldots,j_d}$, we restrict our attention to this subset. By performing this repeatedly, we may assume there is a fixed $d$-tuple $\mathbf{m}$ and a set $\mathcal{B}_{j_1,\ldots,j_d}' \subseteq \mathcal{B}_{j_1,\ldots,j_d}$ such that $\# \mathcal{B}_{j_1,\ldots,j_d}' \geq \# \mathcal{B}_{j_1,\ldots,j_d}/m_1^{\e_1 \delta/2} \cdots m_d^{\e_1 \delta/2}$, $\mathbf{m}$ divides all elements of $\mathcal{B}_{j_1,\ldots,j_d}'$, and there is no $d$-tuple $\mathbf{m}' \neq (1,\ldots,1)$ such that $\mathbf{m}'\mathbf{m}$ divides at least $\# \mathcal{B}_{j_1,\ldots,j_d}'/(m_1' \cdots m_d')^{\e_1 \delta/2}$ elements of $\mathcal{B}_{j_1,\ldots,j_d}'$. Since $m_i \leq 2^{j_i+1} \leq 2XQ^{2/d}$, we have 
    $$\# \mathcal{B}_{j_1,\ldots,j_d}' \geq \frac{\# \mathcal{S}^{\frac{1}{2}-\e_1}}{(2X)^{ \e_1 \delta d/2}Q^{\e_1 \delta}} \geq X^{-\e_1 \delta d/2} \# \mathcal{S}^{\frac{1}{2}-3 \e_1}.$$
    Since $\mathcal{B}_{j_1,\ldots,j_d}' \subseteq \{ \mathbf{b} \in \prod_{i=1}^d [2^{j_i},2^{j_i+1}): m_i \mid b_i \}$ is a set of size $O(2^{j_1+\cdots+j_d}/m_1\cdots m_d)$, we see this also implies $m_1 \cdots m_d < 2^{j_1+\cdots+j_d}/Q^{\delta/4}$. We let
    $$\mathcal{B} := \left\{\mathbf{b}:(m_1b_1,\ldots,m_db_d) \in \mathcal{B}_{j_1,\ldots,j_d}'\right\},$$
    and note $\mathcal{B} \subseteq \prod_{i=1}^d [B_i,2B_i)$, where we set $B_i := 2^{j_i}/m_i$. Let $B = \max_i B_i$, and note that $\# \mathcal{B} \geq (2X)^{-\e_1 \delta d/2}\# \mathcal{S}^{\frac{1}{2}-2\e_1}$ and $B \leq XQ^{2/d}$.\\

     Consider the graph $G=(\mathcal{V},\mathcal{E})$, where the vertex set $\mathcal{V}$ is taken to be $\mathcal{B}$, and the edge set $\mathcal{E}$ is defined by
    $$\mathcal{E} := \left\{ (\mathbf{b}^{(1)},\mathbf{b}^{(2)}) \in \mathcal{B}^2 : \text{exists }1 \leq i \leq d \ \text{ s.t. }\gcd(b_i^{(1)},b_i^{(2)}) \geq B^{\e \delta^2/6r^2} \right\}.$$
    We consider separately two cases.\\

     \textbf{Case 1:} $\# \mathcal{E} \geq \e_1 \# \mathcal{V}^2/2r^2$.\\
    In this case, there are many pairs with a gcd of some size. If we pick a vertex $v$ in $G$ at random, then the expected number of vertices connected to $v$ is at least ${\e_1 \#\mathcal{V}}/{2r^2}$, and so (by the pigeonhole principle) there is some $\mathbf{b}^{(0)} \in \mathcal{B}$ and $1 \leq i \leq d$ such that there are at least ${\e_1 \#\mathcal{B}}/{2r^2d}$ elements $\mathbf{b} \in \mathcal{B}$ with $\gcd(b_i, b_i^{(0)}) > B^{\e \delta^2/6r^2}$. Since there are at most $B^{o(1)}$ divisors of $b_i^{(0)}$, there must be a divisor $m > B^{\e \delta^2/6r^2}$ of $b_i^{(0)}$ such that $m | b_i$ for at least ${\e_1 \#\mathcal{B}}/{2r^2dB^{o(1)}} > {\#\mathcal{B}}/{m^{\e_1 \delta/2}}$ elements $\mathbf{b} \in \mathcal{B}$. But this contradicts the fact that $\mathcal{B}$ is constructed to have no such integers. Thus we must instead have $\#\mathcal{E} < {\e_1 \#\mathcal{V}^2}/{2r^2}$.\\

     \textbf{Case 2:} $\# \mathcal{E}< \e_1 \# \mathcal{V}^2/2r^2$.\\
    Picking $r$ distinct vertices in $G$ uniformly at random, then the expected number of edges between these vertices is less than $\e_1$. In particular, the probability that there are no edges between any of the $r$ chosen vertices is at least $1-\e_1$. Thus, if we define
    $$\mathcal{C}=\left\{ (\mathbf{b}^{(1)},\ldots,\mathbf{b}^{(r)}) \in \mathcal{B}^r: \gcd(b^{(i)}_s,b^{(j)}_s)<B^{\e \delta^2/6r^2} \text{ for } 1 \leq i<j \leq r, \ 1 \leq s \leq d \right\},$$
    then $\# \mathcal{C} \gg_{r,\e} \# \mathcal{B}^r$. Given $(\mathbf{b}^{(1)},\ldots,\mathbf{b}^{(r)}) \in \mathcal{C}$, let
    \begin{align*}
    \mathcal{R}(\mathbf{b}^{(1)},\ldots,\mathbf{b}^{(r)})= &\Big\{ (\mathbf{b}^{(1)'},\ldots,\mathbf{b}^{(r)'}) \in \mathcal{C}: \exists \mathbf{a}^{(1)},\ldots,\mathbf{a}^{(r)},\mathbf{a}^{(1)'},\ldots,\mathbf{a}^{(r)'} \text{ s.t. }\\
    &\gcd(a_i^{(j)},m_ib_i^{(j)})=\gcd(a_i^{(j)'},m_ib_i^{(j)'})=1, \frac{a_i^{(1)}}{b_i^{(1)}}+\cdots+\frac{a_i^{(r)}}{b_i^{(r)}}=\frac{a_i^{(1)'}}{b_i^{(1)'}}+\cdots+\frac{a_i^{(r)'}}{b_i^{(r)'}}, \forall i,j \Big\}.
    \end{align*}
    Abusing notation and defining $\mathcal{R}$ the same as in Lemma \ref{smallgcdfewdenom}, we have $$\# \mathcal{R}(\mathbf{b}^{(1)},\ldots,\mathbf{b}^{(r)}) \leq \prod_{i=1}^d \# \mathcal{R}(b_i^{(1)},\ldots,b_i^{(r)}) \leq  B^{5d\e \delta^2/6}.$$
    For each $\mathbf{b} \in \mathcal{B}$, let $\mathbf{a}(\mathbf{b})$ be an integer tuple where each entry is coprime to that of $\mathbf{m} \mathbf{b}$, and such that $(\mathbf{a}(\mathbf{b}),\mathbf{m} \mathbf{b},\mathbf{h}) \in \mathcal{S}$ for some $\mathbf{b}$. Thus,
    \begin{align*}
        \# \mathcal{A} &\geq \# \left\{ \frac{{a}(\mathbf{b}^{(1)})_i}{m_i b_i^{(1)}}+\cdots+\frac{{a}(\mathbf{b}^{(r)})_i}{m_i b_i^{(r)}}: (\mathbf{b}^{(1)},\ldots,\mathbf{b}^{(r)}) \in \mathcal{C} \right\}\\
        &\geq \sum_{(\mathbf{b}^{(1)},\ldots,\mathbf{b}^{(r)}) \in \mathcal{C}} \frac{1}{\# \mathcal{R}(\mathbf{b}^{(1)},\ldots,\mathbf{b}^{(r)})}\\
        &\geq \frac{\# \mathcal{C}}{B^{5d\e \delta^2/6}} \geq \frac{\# \mathcal{B}^r}{X^{5d\e \delta^2/6}Q^{(5+\frac{\e}{20})\e \delta^2/6}}.
    \end{align*}
    where we recall $r \geq 1>\delta$ and $\# \mathcal{B}>X^{- \e_1 \delta d/2}\# \mathcal{S}^{\frac{1}{2}-3 \e_1} \geq X^{- \e_1 \delta d/2}Q^{(\frac{1}{2}-3 \e_1)\delta}$. This gives condition (2), and we are done.
\end{proof}
 Using this, we obtain an improvement to \cite[][Lemma 7.3]{maynard2021simultaneous}.
\begin{lem}
\label{samedenominators}
    Let $\e>0$, $\delta \in (0,1)$ and $x \in \mathbb{R}_{\geq 1}$. Let $M \in \mathbb{R}^+$ satisfy
    $$M \geq \frac{1}{2}+\frac{\log \e^{-1}}{2 \log 2}.$$
    Let $k,d \geq 2$ be positive integers and let $\alpha_{i,j} \in [0,1)$ for $1 \leq i \leq k$, $1 \leq j \leq d$ be reals. Let $Q$ be large enough in terms of $\delta,d,k,\e$, and let $\e_1,\ldots,\e_k \in (0,1]$ be such that $\Delta=\prod_{i=1}^k \e_i$ satisfies $Q^{(1+\e) \delta} \leq \Delta^{-1} \leq Q^{(2k)^M}$.\\

     Let $\mathcal{S} \subseteq \mathbb{Z}^d \times \mathbb{Z}^d \times \mathbb{Z}^k$ be a set of triples $(\mathbf{a},\mathbf{q},\mathbf{h})$ satisfying:
    \begin{enumerate}
        \item $\gcd(a_j,q_j)=1$ and $q_j \leq x^\e  Q^{1/(d-1)}$ for $j \in \{1,2,\ldots,d\}$.
        \item $h_i \leq \e_i^{-1} \Delta^{-1/(2k)^M}$ for $i \in \{1,2,\ldots,k\}$.
        \item $\# \mathcal{S} \geq Q^{\delta}$.
        \item For each $j \in \{1,2,\ldots,d\}$, we have
        $$\left| h_1 \alpha_{1,j}+\cdots+h_k \alpha_{k,j} - \frac{a_j}{q_j} \right| \leq x^{-\e} \left( \prod_{i=1}^k \e_i \right)^{(2+k/(2k)^M+200\e)/\delta}.$$
    \end{enumerate}
    Then there is a $\mathbf{q}_0 \in \mathbb{Z}^d$ such that at least $\# \mathcal{S}^{1/2}$ of the triples $(\mathbf{a},\mathbf{q},\mathbf{h}) \in \mathcal{S}$ have $\mathbf{q}=\mathbf{q}_0$.
\end{lem}
\begin{proof}
    The proof is similar to that of \cite[][Lemma 7.2]{maynard2021simultaneous}. Choose an integer $r$ such that $Q^{\delta r/(2+k/(2k)^M+10\e)}>x^{\e} \prod_{i=1}^k \e_i^{-1} \geq Q^{\delta r/(2+k/(2k)^M+100 \e)}$. A computation shows $r$ must exist and satisfy $$r \in [(2+k/(2k)^M+10 \e)(1+\e \log_Q x + \e),(2+k/(2k)^M+100 \e)(\e \log_Q x+(2k)^M)/\delta].$$ 
   Therefore we may assume that $Q$ is sufficiently large in terms of $r$.\\
    If $(\mathbf{a}^{(1)},\mathbf{q}^{(1)},\mathbf{h}^{(1)}),\ldots(\mathbf{a}^{(r)},\mathbf{q}^{(r)},\mathbf{h}^{(r)}) \in \mathcal{S}$, then we have for $1\leq s \leq r$ such that for all $j=1,2,\ldots,d$,
    $$\sum_{i=1}^k h_i^{(s)} \alpha_{i,j} = \frac{a_j}{q_j}+O\left( x^{-\e} \left(\prod_{i=1}^k \e_i \right)^{(2+k/(2k)^M +200\e)/\delta} \right)=\frac{a_j^{(s)}}{q_j^{(s)}}+O(Q^{-(1+\e/1000)r}).$$
    Adding these together and observing $rQ^{-\e r/1000} \leq 1$ for $Q$ sufficiently large gives
    $$\frac{a_j^{(1)}}{q_j^{(1)}}+\cdots+\frac{a_j^{(r)}}{q_j^{(r)}}+O(Q^{-r})=\sum_{i=1}^k \alpha_{i,j} \Tilde{h}_i,$$
    where $\Tilde{h}_i=\sum_{s=1}^r h_i^{(s)}$. Since the denominator of $\frac{a_j^{(1)}}{q_j^{(s)}}+\cdots+\frac{a_j^{(r)}}{q_j^{(r)}}$ (when written as a single fraction) is at most $Q^r$, we see that this fraction is uniquely determined by the integers $\Tilde{h}_1,\ldots,\Tilde{h}_k$. Note $|\Tilde{h}_j| \ll r\e_i^{-1} \Delta^{-1/(2k)^M}$, so we have
    \small
    \begin{align*}
    &\# \left\{ \left( \frac{a_1^{(1)}}{q_1^{(1)}}+\cdots+\frac{a_1^{(r)}}{q_1^{(r)}},\ldots,\frac{a_d^{(1)}}{q_d^{(1)}}+\cdots+\frac{a_d^{(r)}}{q_d^{(r)}} \right) : \exists \mathbf{h}^{(1)},\ldots,\mathbf{h}^{(r)} \text{ s.t. } (\mathbf{a}^{(1)},\mathbf{q}^{(1)},\mathbf{h}^{(1)}),\ldots,(\mathbf{a}^{(r)},\mathbf{q}^{(r)},\mathbf{h}^{(r)}) \in \mathcal{S} \right\}\\
    &\leq \# \left\{ (\Tilde{h}_1,\ldots,\Tilde{h}_k) \in \mathbb{Z}^k : |\Tilde{h}_i| \ll r \e_i^{-1} \Delta^{-1/(2k)^M} \text{ for }i=1,2,\ldots,k \right\}\\
    &\leq \left( \prod_{i=1}^k \e_i^{-1} \right)^{1+k/(2k)^M+\e/1000}\\
    &\leq x^{-(1+k/(2k)^M+\e/1000)\e} Q^{\frac{1+k/(2k)^M+\e/1000}{2+k/(2k)^M+100 \e}\delta r} \leq x^{-(1+k/(2k)^M+\e/1000)\e} Q^{\frac{1}{2}-\e}
    \end{align*}
    \normalsize
    Here we used the fact that $r \leq \delta^{-1} \log \prod_{i=1}^k \e_i^{-1}$ and $x^{-\e} Q^{\delta r/(2+k/(2k)^M+100\e)} \leq \prod_{i=1}^k \e^{-1}$, followed by $\prod_{i=1}^k \e_i^{-1} < x^{-\e} Q^{\delta r/(2+k/(2k)^M+10\e)}$.\\

     Note our assumptions satisfy all hypotheses of Lemma \ref{alotorsamedenom}, so by the lemma and also the above bound, (1) must hold, so there must be a $\mathbf{q}_0 \in \mathbb{Z}^d$ such that at least $\# \mathcal{S}^{1/2}$ of the triples $(\mathbf{a},\mathbf{q},\mathbf{h}) \in \mathcal{S}$ have $\mathbf{q}=\mathbf{q}_0$.
\end{proof}
 Therefore, we can prove a proposition akin to \cite[][Proposition 5.2]{maynard2021simultaneous}.
\begin{prop}
\label{5.2}
    Let $\e>0$, and $f_1, \ldots, f_k \in$ $\mathbb{R}[X]$ be polynomials of degree at most $d $ such that $f_1(0)=\cdots=f_k(0)=0$. Put $f_i(X)=\sum_{j=1}^d f_{i, j} X^j$. Let $\e_1, \ldots, \e_k \in(0,1 / 100]$, and put $\Delta=\prod_{i=1}^k \e_i$. Define constants
    $$c_0=1+\e, \quad c_1=1+k/(2k)^M+10\e, \quad c_2=2.5+1/(2k)^{M-1}+500\e,$$
    and let $M \in \mathbb{R}^+$ satisfy
    $$M \geq \frac{1}{2}+\frac{\log \e^{-1}}{2 \log 2}.$$
 Let $Q \leq \Delta^{-c_1d(d-1)}$ be such that there are at least $Q^{1/(1+\e)c_0d(d-1)}$ triples $(\mathbf{a}, \mathbf{q}, \mathbf{h}) \in$ $\mathbb{Z}^d \times \mathbb{Z}^d \times \mathbb{Z}^k$ satisfying:
\begin{enumerate}
\item $q_1=\cdots=q_d$ and $1 \leq q_j \leq  x^\e Q^{1/c_0(d-1)}$ for $1 \leq j \leq d$.
\item $h_i \ll \e_i^{-1} \Delta^{-1 /(2 k)^M}$ for $1 \leq i \leq k$.
\item For each $j \in\{1, \ldots, d\}$ we have
$$
\sum_{i=1}^k h_i f_{i, j}=\frac{a_j}{q_j}+O\left(\frac{Q^{1/d}}{x^{j-\e}}\right) .
$$
\end{enumerate}
Then provided $\Delta^{-c_2d(d-1)}<x$ there is some positive integer $q \leq x^\e Q^{1/c_0(d-1)}$ and at least $Q^{1 / 2(1+\e)c_1d(d-1)}$ pairs $(\mathbf{a}, \mathbf{h}) \in$ $\mathbb{Z}^d \times \mathbb{Z}^k$ such that:
\begin{enumerate}
\item $h_i \ll \e_i^{-1} \Delta^{-2 /(2 k)^M}$ for $i \in\{1, \ldots, k\}$.
\item For each $j \in\{1, \ldots, d\}$ we have
$$
\sum_{i=1}^k h_i f_{i, j}=\frac{a_j}{q}+O\left(\frac{Q^{1/d}}{x^{j-\e}}\right) .
$$
\end{enumerate}
All implied constants depend only on $d,k,\e$ and $M$.
\end{prop}
\begin{proof}
    The proof is nearly identical to \cite[][Proposition 5.2]{maynard2021simultaneous}. By assumption, there is some $Q \leq (\prod_{i=1}^k \e_i^{-1})^{c_1d(d-1)}$ such that there are at least $Q^{1/2c_1d(d-1)}$ triples $(\mathbf{a},\mathbf{q},\mathbf{h}) \in \mathbb{Z}^d \times \mathbb{Z}^d \times \mathbb{Z}^k$ with $\gcd(a_j,q_j)=1$ and $q_j \leq x^\e Q^{1/c_0(d-1)}$ for $j=1,2,\ldots,d$, and with $|h_i| \ll \e_i^{-1} \Delta^{-1/(2k)^M}$ for $i=1,2,\ldots,k$ and with $$\sum_{i=1}^k h_i f_{i,j}=\frac{a_j}{q_j}+O \left( \frac{Q^{1/d}}{x^{j-\e}} \right).$$
    If $Q \leq \Delta^{-1/(2k)^M}$, then we just take one such triple $(\mathbf{a},\mathbf{q},\mathbf{h})$. In this case the triple $(j\mathbf{a},j\mathbf{q},j\mathbf{h})$ for $j=1,2,\ldots,Q$ then give $Q$ relations of the desired type $C_d'>2d(d-1)+1$, since $jh_i\ll Q \e_i^{-1} \Delta^{-1/(2k)^M} \ll \e_i^{-1} \Delta^{-2/(2k)^M}$. Thus we may assume $Q^{(2k)^M} > \Delta^{-1}$.\\

     We now apply Lemma \ref{samedenominators} with $\delta=((1+\e)c_1d(d-1))^{-1}$. Provided $c_2=2.5+1/(2k)^{M-1}+500\e$ and $c_1=1+k/(2k)^{M}+10\e$, we see that the bounds $Q \leq (\prod_{i=1}^k \e_i^{-1})^{c_1d(d-1)}$ and $(\prod_{i=1}^k \e_i^{-1})^{c_2d(d-1)} \leq x$ imply $$\frac{Q^{1/d}}{x^{j-\e}}<x^{-\e}\left(\prod_{i=1}^k \e_i \right)^{c_2(1-2\e)d(d-1)-c_1(d-1)}<x^{-\e}\left(\prod_{i=1}^k \e_i \right)^{(2+k/(2k)^M+200\e)/\delta}$$ for $1 \leq i \leq d$, since
    $$c_2(1-2\e)d(d-1)-c_1(d-1) \geq (2+k/(2k)^M+200\e)(1+\e)c_1d(d-1)$$
    is implied by the fact that
    $$d  \geq \frac{1+k/(2k)^{M}+\e}{(2.5+1/(2k)^{M-1}+500\e)(1-2\e)-(2+k/(2k)^M+200\e)(1+k/(2k)^{M}+10\e)(1+\e)},$$
    and so all hypotheses of Lemma \ref{samedenominators} are satisfied. Therefore, there is some $\mathbf{q}_0 \in \mathbb{Z}^d$ such that at least $Q^{\delta/2}$ of the triples $(\mathbf{a},\mathbf{q},\mathbf{h})$ have $\mathbf{q}=\mathbf{q}_0$. These all give rise to a rational $a_j/q$ due to the condition $q_1=\cdots=q_d$, and so we are done.
\end{proof}
\section{Density Increment} \label{densityincrement}
We also need to alter \cite[][Lemma 8.3]{maynard2021simultaneous} slightly.
\begin{lem}
\label{8.3}
    Let $f_1, \ldots, f_k \in \mathbb{R}[X]$ be polynomials of degree at most $d$ with $f_1(0)=\cdots=$ $f_k(0)=0$. Put $f_i(X)=\sum_{j=1}^d f_{i, j} X^j$.\\

 Let $c>0$, $B_1, \ldots, B_k \geq 1$, $Q \in \mathbb{Z}_{>0}$, and $1 \leq q_0 <x^{\e} Q^{2/d}$. Let $\eta \in[0,1 / 100]$ be such that
$$
\eta<\frac{Q^{1/d}}{x^{1-\e}} .
$$

Let $r \in\{1, \ldots, k\}$ and $\mathbf{h}^{(1)}, \ldots, \mathbf{h}^{(r)} \in \mathbb{Z}^k$ and $\mathbf{a}^{(1)}, \ldots, \mathbf{a}^{(r)} \in \mathbb{Z}^d$ satisfy:
\begin{enumerate}
\item $h_i^{(\ell)} \leq B_i$ for $1 \leq i \leq k$ and $1 \leq \ell \leq r$.
\item For $1 \leq \ell \leq r$ and $1 \leq j \leq d$ we have
$$
\left|\sum_{i=1}^k h_i^{(\ell)} f_{i, j}-\frac{a_j^{(\ell)}}{q_0}\right| \leq \eta^j .
$$
\item Put $\tilde{h}_i^{(\ell)}=h_i^{(\ell)} / B_i$ for $1 \leq \ell \leq r$ and $1 \leq i \leq k$. We have
$$
\left\|\tilde{\mathbf{h}}^{(1)} \wedge \cdots \wedge \tilde{\mathbf{h}}^{(r)}\right\| \asymp\left\|\tilde{\mathbf{h}}^{(1)}\right\|_{\infty} \cdots\left\|\tilde{\mathbf{h}}^{(r)}\right\|_{\infty}
$$
and
$$
\left\|\tilde{\mathbf{h}}^{(1)}\right\|_{\infty} \cdots\left\|\tilde{\mathbf{h}}^{(r)}\right\|_{\infty} \ll \frac{1}{Q^{1 / cd(d^2-1)}} .
$$
\end{enumerate}
Then there is an integer $k^{\prime}<k$, real polynomials $g_1, \ldots, g_{k^{\prime}} \in \mathbb{R}[X]$ of degree at most $d$ with $g_1(0)=\cdots=g_k(0)=0$ and quantities $B_1^{\prime}, \ldots, B_{k^{\prime}}^{\prime} \geq 2$ and $y<x$ such that:
\begin{enumerate}
\item If there is an integer $n^{\prime}<y$ such that
$$
\left\|g_i\left(n^{\prime}\right)\right\|<\frac{1}{B_i^{\prime}} \quad \text { for all } 1 \leq i \leq k^{\prime}
$$
then there is an integer $n<x$ such that
$$
\left\|f_i(n)\right\|<\frac{1}{B_i} \quad \text { for all } 1 \leq i \leq k .
$$
\item We have
$$
\frac{y}{\left(B_1^{\prime} \cdots B_{k^{\prime}}^{\prime}\right)^{d(d-1)(4.5c+1.5-1/k'^3)}} \gg \frac{x^{1-\e}}{\left(B_1 \cdots B_k\right)^{d(d-1)(4.5c+1.5-1/k^4-1/k^4)}}.
$$
\end{enumerate}
\end{lem}
\begin{proof}
    The majority of the proof is identical to the corresponding proof in \cite{maynard2021simultaneous}. Due to the proof being quite long, we only include the parts where it differs. The following is a modified version of the proof, bottom half of Page 24 of \cite{maynard2021simultaneous} onwards.\\

    We see that $g_i$ are polynomials of degree at most $d$ with $g_i(0)=0$ since $\tilde{f}_i$ are. Finally, we put $y=\delta x^{1-\e} \min_i \lVert \tilde{\mathbf{h}}_i \rVert_\infty/q_0Q^{1/d} D_2$, and note that since $\eta<Q^{1/d}/x^{1-\e}$, we have $y<\delta \min_i \lVert \tilde{\mathbf{h}}_i \rVert_\infty/\eta q_0D_2$. Putting everything together, we see that there is an $n'<y$ such that
    $$\lVert g_i(n') \rVert \leq \frac{1}{B_i'}$$
    for $1 \leq i \leq k'=k-r$, then there is an $n=n'q_0D_2$ with $n<x$ and $n<\delta \min_i \lVert \tilde{\mathbf{h}}_i \rVert_\infty/\eta$ such that
    $$\lVert f_i(n) \rVert \leq \frac{1}{B_i}$$
    for $1 \leq i \leq k$. Thus we are left to verify the size estimates with this choice of $B_1',\ldots,B_{k-r}'$ and $y$. We have
    $$\prod_{i=1}^{k-r} B_i' = \frac{\lVert \mathbf{z}_1 \rVert_\infty \cdots \lVert \mathbf{z}_{k-r} \rVert_\infty \prod_{i=r+1}^k B_i}{\delta^{2(k-r)}} \ll \frac{D_1 \prod_{i=r+1}^k B_i}{D_2}.$$
    This implies that (with $C_2$ a constant depending only on $d$ chosen later)
    $$\frac{y}{(\prod_{i=1}^{k-r} B_i')^{C_2}}=\frac{\delta x^{1-\e} \min_i \lVert \tilde{\mathbf{h}}_i \rVert_\infty}{q_0Q^{1/d}D_2 (\prod_{i=1}^{k-r} B_i')^{C_2}} \gg \frac{x^{1-2\e} D_2^{C_2-1} \min_i \lVert \tilde{\mathbf{h}}_i \rVert_\infty}{Q^{3/d}D_1^{C_2}(\prod_{i=r+1}^{k} B_i)^{C_2}}.$$
    Recall that $D_1=\det(H_1) \asymp \lVert \mathbf{h}_1 \rVert_\infty \cdots \lVert \mathbf{h}_r \rVert_\infty \ll B_1 \cdots B_r/Q^{1/cd(d-1)^2}$, and that $\min_i \lVert \tilde{\mathbf{h}}_i \rVert_\infty \gg \prod_{i=1}^r \lVert \tilde{\mathbf{h}}_i \rVert_\infty \gg D_1/(B_1 \cdots B_r)$. This gives
    \begin{align*}
        \frac{y}{(\prod_{i=1}^{k-r} B_i')^{C_2}} &\gg \frac{x^{1-2\e}}{(\prod_{i=r+1}^{k} B_i)^{C_2}} \frac{D_2^{C_2-1}}{D_1^{C_2-1} Q^{3/d} B_1 \cdots B_r}\\
        &\gg \frac{x^{1-2\e}}{(\prod_{i=1}^{k} B_i)^{C_2}} \cdot Q^{(C_2-1)/cd(d-1)^2-3/d} D_2^{C_2-1}.
    \end{align*}
    Finally, we choose $C_2=d(d-1)(4.5c+1.5-1/(k-r)^3)$. Since $D_2,Q \geq 1$ and $$(C_2-1)/cd(d-1)^2-3/d \geq 0,$$
   which gives
    $$\frac{y}{(\prod_{i=1}^{k-r} B_i')^{d(d-1)(4.5c+1.5-1/(k-r)^3)}} \gg \frac{x^{1-\e}}{(\prod_{i=1}^r B_i)^{d(d-1)(4.5c+1.5-1/(k-r)^3)}}.$$
    Since $k>r \geq 1$, we have $1/(k-r)^3 \geq 1/k^3+1/k^4$, and so we are done.
\end{proof}
There are minimal changes to the proof of \cite{maynard2021simultaneous} from here onwards. For completeness, we include the remaining statements and proofs.
\begin{lem} \label{8.1}
    Let $\eta>0$ be sufficiently small in terms of $k$ and $d$. Let $B_1 ,\ldots,B_k>1$ satisfy $\prod_{i=1}^k B_i \leq \eta^{-1/2}$ and $\beta_{i,j} \in \mathbb{R}$ for $1 \leq i \leq k$, $1 \leq j \leq d$. Let $\mathcal{R}$ be the region in $\mathbb{R}^{k+d}$ defined by
        \begin{align*}
            \mathcal{R}=\{ &(h_1,\ldots,h_k,a_1,\ldots,a_d) \in \mathbb{R}^{k+d}:\\
            &\left| \sum_{i=1}^k h_i \beta_{i,j}-a_j \right| \leq \eta^j \enspace \forall 1 \leq j \leq d, \, |h_i| \leq B_i \enspace \forall 1 \leq i \leq k\},
        \end{align*}
        and assume that $\#(\mathcal{R} \cap \mathbb{Z}^{k+d})=N$ is sufficiently large in terms of $k$ and $d$. Then $\exists 1 \leq r \leq k$ and $\mathbf{h}^{(1)},\ldots,\mathbf{h}^{(r)} \in \mathbb{Z}^k$ and $\mathbf{a}^{(1)},\ldots,\mathbf{a}^{(r)} \in \mathbb{Z}^d$ such that:
        \begin{enumerate}
            \item $\forall j \in \{1,\ldots,r\}$, $(h_1^{(j)},\ldots,h_k^{(j)},a_1^{(j)},\ldots,a_d^{(j)}) \in \mathcal{R} \cap \mathbb{Z}^{k+d}.$
            \item $\lVert \mathbf{h}^{(1)} \wedge \cdots \wedge \mathbf{h}^{(r)} \rVert \asymp_{k,d} \lVert \mathbf{h}^{(1)} \rVert_\infty \cdots \lVert \mathbf{h}^{(r)} \rVert_\infty$
            \item $\lVert \mathbf{h}^{(1)} \rVert_\infty \cdots \lVert \mathbf{h}^{(r)} \rVert_\infty \ll_{k,d} {B_1 \cdots B_k}/{N^{1/(d+1)}}.$
        \end{enumerate}
        All implied constants depend at most on $k$ and $d$.
\end{lem}
\begin{proof}
    This is \cite[][Lemma 8.1]{maynard2021simultaneous}.
\end{proof}
\begin{prop}
\label{5.3}
    Let $f_1,\ldots,f_k \in \mathbb{R}[X]$ be polynomials of degree at most $d \geq 2$ such that $f_1(0)=\cdots=f_k(0)=0$. Put $f_i(X)=\sum_{j=1}^d f_{i,j} X^j$. Let $\e_1,\ldots,\e_k \in (0,1/100]$, and put $\Delta=\prod_{i=1}^k \e_i$. For any $\e>0,M \geq 4$, define constants
    $$c_0=1+\e, \quad c_1=1+k/(2k)^M+10\e, \quad c_2=2.5+1/(2k)^{M-1}+500\e.$$
    Suppose $\Delta^{-c_2d(d-1)} \leq x$, and let $Q \leq \Delta^{-c_1d(d-1)}$. Let $q$ be a positive integer with $q<x^\e Q^{1/c_0(d-1)}$. Let $\mathcal{S}$ be the set of pairs $(\mathbf{a},\mathbf{h}) \in \mathbb{Z}^d \times \mathbb{Z}^k$ such that for $j \in \{1,\ldots,d\}$, we have
    $$\left| \sum_{i=1}^k h_i f_{i,j} -\frac{a_j}{q} \right| \ll \frac{Q^{1/d}}{x^{j-\e}},$$
    and such that $|h_i| \ll \e_i^{-1} \Delta^{-2/(2k)^M}$. Assume that $\# \mathcal{S} > Q^{1/2(1+\e)c_1d(d-1)}$.\\

    Then there is an integer $k'<k$, polynomials $g_1,\ldots,g_{k'} \in \mathbb{R}[X]$ of degree at most $d \geq 2$ with $g_1(0)=\cdots=g_{k'}(0)=0$ and quantities $\e_1',\ldots,\e_{k'}' \in (0,1/100]$ and $y<x$ such that for $c=9(1+\e)c_1+1.5$:
    \begin{enumerate}
        \item If there is an integer $n'<y$ such that for all $1 \leq i \leq k'$
        $$\lVert g_i(n') \rVert<\e_i' $$
        then there is an integer $n<x$ such that for all $1 \leq i \leq k$,
        $$\lVert f_i(n) \rVert<\e_i.$$
        \item We have
        $$y(\e_1' \cdots \e_{k'}')^{d(d-1)(c-1/k'^3)} \gg x^{1-\e}(\e_1 \cdots \e_k)^{d(d-1)(c-1/k^3)}.$$
    \end{enumerate}
    All implied constants depend only on $k$, $d$ and $\e$.
\end{prop}
\begin{proof}
    Taking $\beta_{i,j}=f_{i,j}$, $\eta=Q^{1/d}/x$ and $B_i \ll \e_i^{-1} \Delta^{-2/(2k)^M}$, Lemma \ref{8.1} shows that we can find a subset of essentially orthogonal generators, and we can apply Lemma \ref{8.3}. It suffices to verify the density estimates. Indeed,
    \begin{align*}
        y(\e'_1 \cdots \e'_{k'})^{d(d-1)(c-1/k'^3)} &= \frac{y}{(B_1' \cdots B_{k'}')^{d(d-1)(c-1/k'^3)}}\\
        &\gg \frac{x}{(B_1 \cdots B_k)^{d(d-1)(c-1/k^4-1/k^4)}}\\
        &= x(\e_1 \cdots \e_k\Delta^{2/(2k)^M})^{d(d-1)(c-1/k^4-1/k^4)}\\
        &\gg x(\e_1 \cdots \e_k)^{d(d-1)(c-1/k^3)}.
    \end{align*}
    Note here $M \geq 4$ is required for the last step.
\end{proof}
\section{Putting Everything Together} \label{puttingeverythingtogether}
As in \cite{maynard2021simultaneous}, combining Propositions \ref{5.1}, \ref{5.2}, \ref{5.3}, we obtain the following version of \cite[][Proposition 5.4]{maynard2021simultaneous}.
\begin{prop}
\label{5.4}
    Let $\e>0$, and $d,k$ be positive integers. Define constants
    $$c_0=1+\e, \quad  c_2=2.5+1/(2k)^{M-1}+100\e,$$
    and let $M \in \mathbb{R}^+$ satisfy
    $$M \geq \max \left\{4, \frac{1}{2}+\frac{\log \e^{-1}}{2 \log 2} \right\}.$$
    Let $c_2=2.5+1/(2k)^{M-1}+100\e$, and $C_{d,k}$ be a constant depending only on $d$, $k$ (and implicitly on $\e,M$) such that the following holds.\\

    Let $f_1,\ldots,f_k \in \mathbb{R}[X]$ be polynomials of degree at most $d$ such that $f_1(0)= \cdots = f_k(0)=0$. Put $f_i(X)=\sum_{j=1}^d f_{i,j} X^j$. Let $\e_1,\ldots,\e_k \in (0,1/100]$, and put $\Delta=\prod_{i=1}^k \e_i$. Let $\Delta^{-1} \leq x^{1/c_2d(d-1)}$.\\

    If there is no positive integer $n<x$ such that for all $i=1,2,\ldots,k$,
    $$\lVert f_i(n) \rVert<\e_i,$$
    then there is a positive integer $k'<k$ and polynomials $g_1,\ldots,g_{k'} \in \mathbb{R}[X]$ of degree at most $d$, with $g_1(0)=\cdots=g_{k'}(0)=0$ and reals $\e_1',\ldots,\e_{k'}' \in (0,1/100]$ and $y \in \mathbb{R}$ with $y<x$ such that both of the following hold:
    \begin{enumerate}
        \item There is no positive integer $n'<y$ such that for all $i=1,2,\ldots,k'$,
        $$\lVert g_i(n') \rVert<\e_i'.$$
        \item We have
        $$y(\e_1' \cdots \e_{k'}')^{d(d-1)(c-1/k'^3)} \geq \frac{1}{C_{d,k}} x^{1-\e}(\e_1 \cdots \e_k)^{d(d-1)(c-1/k^3)},$$
        where $c=9(1+\e)(1+k/(2k)^M+10\e)+1.5$.
    \end{enumerate}
    All implied constants depend only on $k$, $d$ and $\e$.
\end{prop}
\begin{thm} \label{premaintheorem}
Let $k, d \in \mathbb{Z}^+,\, \e>0,$ and $M \in \mathbb{R}^+$ satisfy
$$M \geq \max \left\{4, \frac{1}{2}+\frac{\log \e^{-1}}{2 \log 2} \right\}.$$
There is a constant $C_{d, k}>2$ depending only on $d$, $k$ (and $\e,M$) such that the following holds.

Let $f_1, \ldots, f_k \in \mathbb{R}[X]$ be polynomials of degree at most $d$ such that $f_1(0)=\cdots=$ $f_k(0)=0$. Let $\e_1, \ldots, \e_k \in(0,1 / 100]$, and put $\Delta=\prod_{i=1}^k \e_i$. Define $$c_2=10.5+9k/(2k)^M+\e.$$

If $\Delta^{-1} \leq x^{1 / c_2d(d-1)}$ and $x>C_{d, k}$, then there is a positive integer $n<x$ such that
$$
\left\|f_i(n)\right\| \leq \e_i \quad \text { for all } i \in\{1, \ldots, k\} \text {. }
$$
\end{thm}
\begin{proof}
    The proof is essentially identical to \cite[][Theorem 1.1]{maynard2021simultaneous}. Assume for a contradiction that there is no positive $n<x$ such that $\lVert f_i(n) \rVert \leq \e_i$ for all $i=1,2,\ldots,k$. Let $C_{d,k}$ be the constant in Proposition \ref{5.4}. Let $C_0=\sup_{j \leq k} C_{d,j}$.\\
    
    Define a \emph{System} to be a tuple $(k,\mathbf{g},\bm{\delta},y)$ consisting of
    \begin{enumerate}
        \item A positive integer $k$.
        \item A $k$-tuple $\mathbf{g}$ of real polynomials $(g_1,\ldots,g_k)$ of degree at most $d$ satisfying $g_1(0)=\cdots=g_k(0)=0$.
        \item A $k$-tuple $\bm{\delta}$ of reals $(\delta_1,\ldots,\delta_k)$ with $\delta_i \in (0,1/100]$ for all $i \in \{1,2,\ldots,k\}$.
        \item A real $y$ such that there is no positive integer $n<y$ satisfying
        $$\lVert g_i(n) \rVert \leq \delta_i$$
        for all $i=1,2,\ldots,k$.
    \end{enumerate}
    Let $c=9(1+\e)(1+k/(2k)^M+\e)+1.5$. Given a system $(k,\mathbf{g},\bm{\delta},y)$, let $\Delta(\bm{\delta})=\prod_{i=1}^k \delta_i$. By Proposition \ref{5.4}, if a system $(k_j,\mathbf{g}_j,\bm{\delta}_j,y_j)$ satisfies $\Delta(\bm{\delta}_j)^{-c_2d(d-1)}<y_j$, then there is a system $(k_{j+1},\mathbf{g}_{j+1},\bm{\delta}_{j+1},y_{j+1})$ such that $k_{j+1}<k_j$, $y_{j+1} \leq y_j$,
    $$y_{j+1} \Delta(\bm{\delta}_{j+1})^{d(d-1)(c-1/k_{j+1}^2)} \geq \frac{y_j^{1-\e} \Delta(\bm{\delta}_j)^{d(d-1)(c-1/k_{j}^2)}}{C_0},$$
    Note that if
    \begin{equation} \label{betterinductionhypothesis}
    y_j^{1-k_j\e} \Delta(\bm{\delta}_j)^{d(d-1)(c-1/k_{j}^2)} \geq C_0^{k_j},
    \end{equation}
    then we have
    $$y_{j+1}^{1-k_{j+1}\e} \Delta(\bm{\delta}_{j+1})^{d(d-1)(c-1/k_{j+1}^2)} \geq y_{j+1}^{-k_{j+1}\e} \frac{y_j^{1-\e} \Delta(\bm{\delta}_j)^{d(d-1)(c-1/k_{j}^2)}}{C_0} \geq C_0^{k_j-1} \geq C_0^{k_{j+1}}.$$
    Therefore, if (\ref{betterinductionhypothesis}) holds, we can find infinite systems with $k_j$ a strictly descending list of positive integers, which is impossible. Therefore, there does not exist a system satisfying (\ref{betterinductionhypothesis}).\\

    Therefore, given an integer $k$, polynomials $\mathbf{f}=(f_1,\ldots,f_k) \in \mathbb{R}[X]^k$ of degree at most $d$ satisfying $f_1(0)=\cdots=f_k(0)=0$, quantities $\bm{\e}=(\e_1,\ldots,\e_k) \in (0,1/100]^k$ and a real $x$ satisfying $x^{1-k\e} \Delta^{d(d-1)(c-1/k^2)} \geq C_0^{k}$, $(k,\mathbf{f},\bm{\e},x)$ cannot form a system, so there must exist some $n<x$ such that
    $\lVert f_i(n) \rVert \leq \e_i$ for all $i=1,2,\ldots,k$.
\end{proof}
 By taking $\e_1=\cdots=\e_k=x^{-1/c_2d(d-1)k}$, and for any $\e>0$ letting $M$ large such that $2^M \geq \e^{-1}$, we arrive at the following:
\begin{cor} \label{maincor}
    Let $f_1,\ldots,f_k \in \mathbb{R}[X]$ be polynomials of degree at most $d$ such that $f_1(0)=\cdots=f_k(0)=0$. Then there is a positive integer $n<x$ such that for all $\e>0$,
    $$\lVert f_i(n) \rVert \ll_{d,k,\e} x^{-1/(10.5+\e)kd(d-1)}$$
    for all $i=1,2,\ldots,k$.
\end{cor}
\newpage
\nocite{*}
\setcitestyle{numbers}
\bibliographystyle{plainnat}
\bibliography{bibliography.bib}
\end{document}